\documentclass[10pt,reqno]{amsart}
\usepackage{amsmath}
\usepackage{amsthm}
\usepackage{amssymb}
\usepackage{amsfonts}
\usepackage{latexsym}
\usepackage{graphicx}
\usepackage{color}

\usepackage{hyperref} 
\usepackage{times} 
\usepackage{fullpage} 
\usepackage{multicol} 
\usepackage{nonfloat}
\usepackage{lipsum} 
\newcommand{\be}{\begin{equation}}
\newcommand{\bea}{\left[ \begin{array}}
\newcommand{\ee}{\end{equation}}
\newcommand{\eea}{\end{array} \right]}
\newcommand{\ba}{\begin{array}}
\newcommand{\ea}{\end{array}}

\newtheorem{theorem}{Theorem}
\newtheorem{corollary}{Corollary}

\newcommand{\true}{\rm true}
\newcommand{\reg}{\rm reg}
\newcommand{\real}{\mathbb R}

\newcommand{\bSigma}{{\boldsymbol{\Sigma}}}

\newcommand{\bA}{{\bf A}}
\newcommand{\bB}{{\bf B}}
\newcommand{\bC}{{\bf C}}

\newcommand{\bI}{{\bf I}}

\newcommand{\bL}{{\bf L}}
\newcommand{\bM}{{\bf M}}

\newcommand{\bP}{{\bf P}}

\newcommand{\bR}{{\bf R}}

\newcommand{\bT}{{\bf T}}
\newcommand{\bU}{{\bf U}}
\newcommand{\bV}{{\bf V}}
\newcommand{\bW}{{\bf W}}
\newcommand{\bX}{{\bf X}}

\newcommand{\bee}{{\bf e}}
\newcommand{\baa}{{\bf a}}
\newcommand{\bb}{{\bf b}}

\newcommand{\bd}{{\bf d}}

\newcommand{\bt}{{\bf t}}
\newcommand{\bu}{{\bf u}}
\newcommand{\bv}{{\bf v}}

\newcommand{\bx}{{\bf x}}

\DeclareMathOperator*{\argmin}{argmin}

\def\I0{\bigl[ \, I_m \, , \, 0_{m\times p} \, \bigr]}
\def\0I{\bigl[ \, 0_{n\times p} \, I_p \, \bigr]}

\newcommand{\xvec}{{\mathbf x}}
\newcommand{\bvec}{{\mathbf b}}
\newcommand{\evec}{{\mathbf e}}
\newcommand{\vvec}{{\mathbf v}}
\newcommand{\uvec}{{\mathbf u}}

\allowdisplaybreaks
\begin{document}
	\title{The Image Deblurring Problem: Matrices, Wavelets, and Multilevel Methods}

	\author{David Austin$^{1}$, Malena I. Espa\~nol${^2}$, and Mirjeta Pasha$^{2}$}
	
	\address{$^1$ Department of Mathematics, Grand Valley State University, Allendale, MI.}
	\address{
	$^2$ School of Mathematical and Statistical Sciences, Arizona State University, Tempe, AZ.}
	\email{austind@gvsu.edu, malena.espanol@asu.edu, mpasha@asu.edu}
	
\begin{abstract}
The image deblurring problem consists of reconstructing images from blur and noise contaminated available data. In this AMS Notices article, we provide an overview of some well known numerical linear algebra techniques that are use for solving this problem. In particular, we start by carefully describing how to represent images, the process of blurring an image and modeling different kind of added noise. Then, we present regularization methods such as Tikhonov (on the standard and general form), Total Variation and other variations with sparse and edge preserving properties. Additionally, we briefly overview some of the main matrix structures for the blurring operator and finalize presenting multilevel methods that preserve such structures. Numerical examples are used to illustrate the techniques described.
\end{abstract}
	
\maketitle

\begin{multicols}{2}
\section{Introduction}

After the launch of the Hubble Space Telescope in 1990, astronomers were gravely disappointed by the quality of the images as they began to arrive.  Due to miscalibrated testing equipment, the telescope's primary mirror had been ground to a shape that differed slightly from the one intended resulting in the misdirection of incoming light as it moved through the optical system.  The blurry images did little to justify the telescope's \$1.5 billion price tag.

Three years later, space shuttle astronauts installed a specially designed corrective optical system that essentially fixed the problem and yielded spectacular images (see Figure \ref{fig:hubble}).  In the meantime, mathematicians devised several ways to convert the blurry images into high-quality images.  The process of mathematical deblurring is the focus of this article.

\begin{myfigure}{
\centering
    \includegraphics{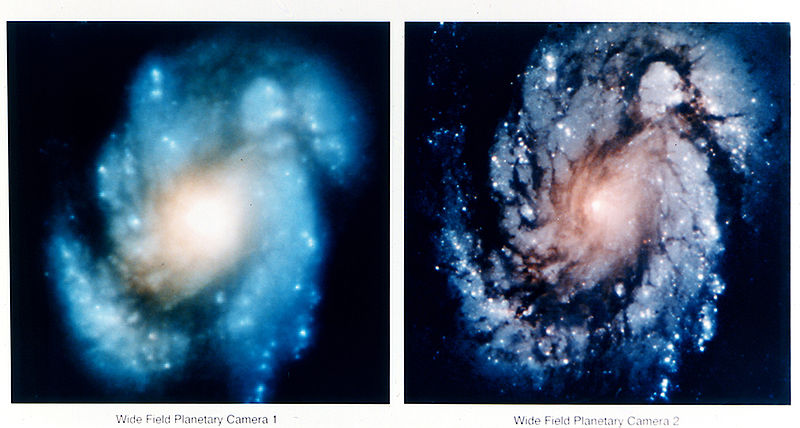}}
\end{myfigure}
\figcaption{Hubble's view of the M100 galaxy, soon after launch on the left and after corrective optics were installed in 1993.  NASA / ESA \label{fig:hubble}
}
Many factors can cause an image to become blurred, such as motion of the imaging device or the target object, errors in focusing, or the presence of atmospheric turbulence \cite{ellappan2017reconstruction}.
Indeed, the need for image deblurring goes beyond Hubble's story. For instance, image deblurring is widely used in many applications, such as pattern recognition, computer vision, and machine intelligence. Moreover, image deblurring shares the same mathematical formulation as other imaging modalities. For instance, in many cases we have only limited opportunities to capture an image; this is particularly true of medical images, such as computerized tomography (CT), proton computed tomography (pCT), and magnetic resonance imaging (MRI), for which equipment and patient availability are scarce resources.  In cases such as these, we need a way to extract meaningful information from noisy images that have been imperfectly collected.

In this article, we will describe a mathematical model of how digital images become blurred as well as several mathematical issues that arise when we try to undo the blurring.  While blurring may be effectively modeled by a linear process, we will see that deblurring is not as simple as inverting that linear process.  Indeed, deblurring belongs to an important class of problems known as {\em discrete ill-posed problems}, and we will introduce some techniques that have become standard for solving them.  In addition, we will describe some additional structures in the linear operators that cause the required computations to be feasible.
\section{Digital images and blurring}

We start by describing digital images and a process by which they become blurred. 

As illustrated in Figure \ref{fig:image-creation}, 
the lens of a digital camera directs photons entering the camera onto a charge-coupled device (CCD), which consists of a rectangular $p\times q$ array of detectors. Each detector in the CCD converts a count of the photons into an electrical signal that is digitized by an analog-to-digital converter (ADC).  The result is a digital image stored as a $p\times q$ matrix whose entries represent the intensities of light recorded by each of the CCD's detectors. 
\begin{myfigure}{
  \centering
    \includegraphics[scale=0.55]{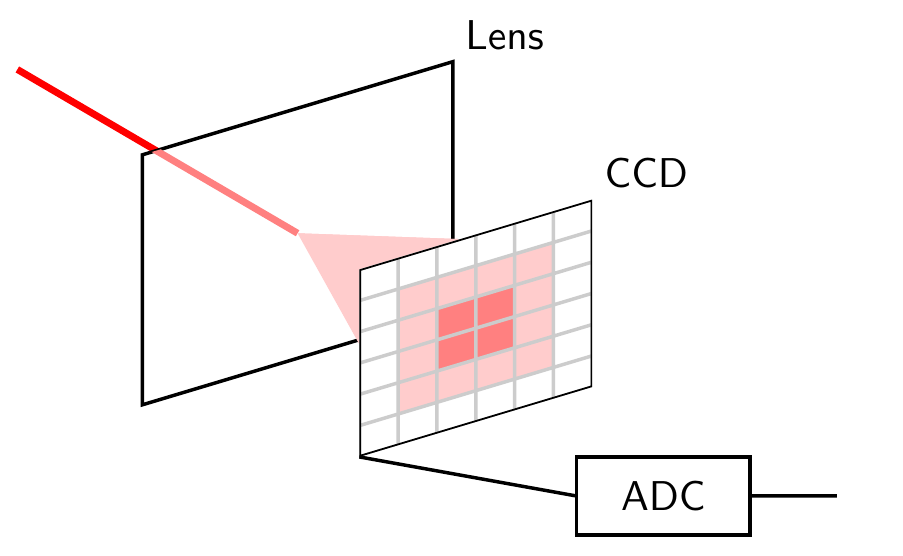}
  }
\end{myfigure}
\figcaption{
  A simple model of how a digital image is created.
  \label{fig:image-creation}
}
A grayscale image is represented by a single matrix with integer entries describing the brightness at each location. A color image is represented by three matrices that describe the colors in terms of their red, green, and blue constituents. Of course, we may see these matrices by zooming in on a digital image until we see individual pixels.

Perhaps due to imperfections in the camera's lens or the lens being improperly focused,
it is inevitable that photons intended for one pixel bleed over into adjacent pixels, and this leads to blurring of the image.  
To illustrate, we will consider grayscale images comprised of arrays of
$64\times 64$ pixels.  In Figure \ref{fig:single-pixel}, the image on
the left shows a single pixel illuminated while on the right we see
how photons intended for this single pixel have spilled over into
adjacent pixels to create a blurred image.

\begin{myfigure}{
  \centering
    \includegraphics[scale=0.45]{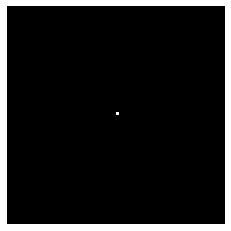}
    \includegraphics[scale=0.45]{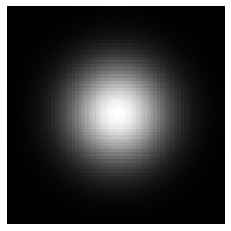}
  }
\end{myfigure}
\figcaption{
  The intensity from a single pixel, shown on the left, is
  spread out across adjacent pixels according to a Gaussian blur, as seen
  on the right.
  \label{fig:single-pixel}
}

There are several models used to describe blurring, but a simple one
that we choose here has the light intensity contained in a single
pixel $(i, j)$ spilling over into an adjacent pixel $(k, l)$ according
to the Gaussian
\begin{equation}
    \label{eq:Gaussian-blur}
\begin{aligned}
  & \frac1N\exp\left(-\frac12\left(\frac{k-i}{s}\right)^2 -
  \frac12\left(\frac{l-j}{s}\right)^2\right) \\
& =
\frac1N\exp\left(-\frac12\left(\frac{k-i}{s}\right)^2\right)
\exp\left(-\frac12\left(\frac{j-k}{s}\right)^2\right),
\end{aligned}
\end{equation}
where $s$ is a parameter that controls the spread in the intensity
and $N$ is a normalization constant so that the total intensity
sums to 1.  As we will see later, the fact that the two-dimensional
Gaussian can be written as a product of one-dimensional Gaussians has important
consequences for our ability to efficiently represent the blurring
process as a linear operator.

Though we visually experience a grayscale image 
as a $p\times q$ matrix of pixels $\bX$, we will mathematically represent an image
as a $pq$-dimensional vector $\xvec$ by stacking the columns of $\bX$ on top of one another.
That is, 
$\xvec = {\rm vec}(\bX)
= (\bX_{11},\ldots,\bX_{p1},\bX_{21},\ldots,\bX_{p2},\ldots,
\bX_{1q},\ldots,\bX_{pq})^T$ with $\bX_{ij}$ being the intensity value of the pixel at row $i$ and column $j$.
The blurring process is linear as the number of photons that arrive at one pixel is the sum of the number of misdirected photons intended for nearby pixels.  Consequently, there is a blurring matrix $\bA$ that blurs the image $\xvec$ into the image $\bvec = \bA\xvec$.  {\em Deblurring} refers to the inverse process of recovering the original image $\xvec$ from its blurred counterpart $\bvec$. 

Each column of $\bA$ is the result of blurring a single pixel, which means that each column of $\bA$ is formed from the image on the right of Figure \ref{fig:single-pixel} by translating to a different center. When that center is near the edge of the image, some photons will necessarily be lost outside the image, and there are a few options for how to incorporate this fact into our model. In real-life settings, it is possible to have knowledge only over a finite region, the so called Field of View (FOV), that defines the range that a user can see from an object. It is then necessary to make an assumption on what is outside the FOV by means of the boundary conditions. For instance, one option to overcome the loss outside the FOV is to simply accept
that loss, in which case we say that the matrix has zero boundary conditions. This has the effect of assuming that the image is black (pixel values are zero) outside the FOV, which can lead to an artificial black border around a deblurred image (see the  image on the left on Figure \ref{fig: boundary}).

In some contexts, it can be advantageous to assume reflexive boundary conditions, which assumes that the photons are reflected back onto the image. In other scenarios of interest, periodic boundary conditions,  which assumes the lost photons reappear on the opposite side of the image as if the image repeats itself indefinitely in all directions outside the FOV, are a suitable fit. Nevertheless, in practical settings, we periodically extend only some pixel values close to the boundary (see the image on the right on Figure \ref{fig: boundary}).
\begin{myfigure}{
    \centering
    \includegraphics[scale=0.45]{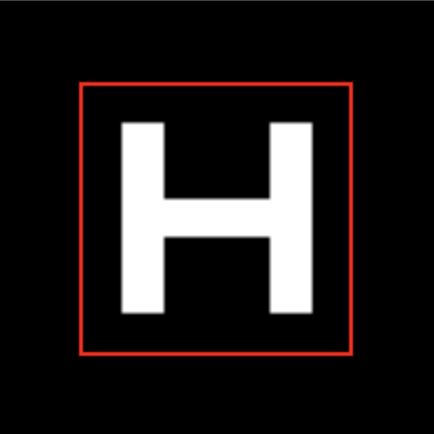}
    \includegraphics[scale=0.45]{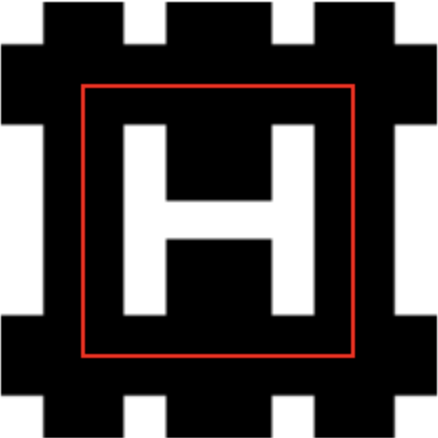}
  }
\end{myfigure}
\figcaption{ Image with assumed zero boundary conditions is shown on the left and one with assumed periodic boundary condition is shown on the right. The red box represents the FOV. 
 \label{fig: boundary}
}
\section{Adding noise}
Let us consider the grayscale image $\xvec^{\true}$ shown on the left of Figure \ref{fig:x-true} and its blurred version $\bvec^{\true} = \bA\xvec^{\true}$ shown on the right.  

\begin{myfigure}{
    \centering
    \includegraphics[scale=0.45]{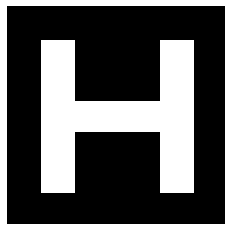}
    \includegraphics[scale=0.45]{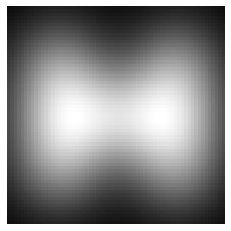}
  }
\end{myfigure}
\figcaption{
  An image $\bx^{\true}$ on the left is blurred to obtain $\bb^{\true}$ on
  the right.
  \label{fig:x-true}
}

If we had access to $\bvec^{\true}$, it would be easy enough to recover $\xvec^{\true}$ by simply solving the linear system $\bA\xvec = \bvec^{\true}$. However, the conversion of photon counts into an electrical signal by the CCD and then a digital reading by the ADC introduces electrical noise into the image that is ultimately recorded.  The recorded image $\bvec$ is therefore a noisy approximation to the true blurred image $\bvec^{\true}$, so we write $\bb = \bb^{\rm true}+\bee$, where  $\bee$ is the noise vector. There are various models used to describe the kind of noise added.  For instance, we might assume that the noise is Gaussian white noise, which means that the entries in $\evec$ are sampled from a normal distribution with mean zero. In our example, we assume that $\|\evec\|_2 = 0.001\|\bvec^{\true}\|_2$, that is, the white noise is about 0.1\% of the image $\bvec^{\true}$.  As seen in Figure \ref{fig:btrue-v-b}, this level of noise cannot be easily detected. 

\begin{myfigure}{
    \centering
    \includegraphics[scale=0.45]{Figures/64x64Case/b_true.png}
    \includegraphics[scale=0.45]{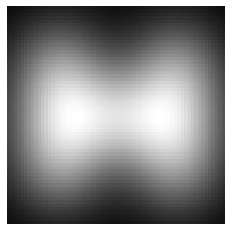}
  }
\end{myfigure}
\figcaption{
  The blurred image $\bvec^{\true}$ on the left with a small amount of Gaussian white noise added to obtain $\bvec$ on the right.
  \label{fig:btrue-v-b}
}

Two other models of noise are illustrated in Figure \ref{fig:other-noise}.  Under low light intensities, as encountered in astrophotography, the number of photons that arrive on the CCD while the image is exposed may differ from the number expected.  A Poisson distribution provides an effective description of the resulting image as seen on the left. If the digital image is transmitted over a communication channel, some of the transmitted bits may be corrupted in transmission resulting in ``salt and pepper" noise demonstrated on the right.

\begin{myfigure}{
    \centering
    \includegraphics[scale=0.45]{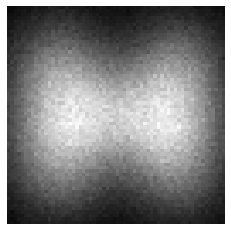}
    \includegraphics[scale=0.45]{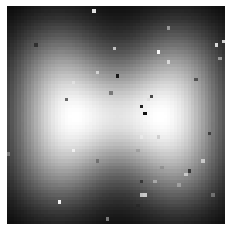}
  }
\end{myfigure}
\figcaption{
  Poisson noise is added to the blurred image on the left and salt and pepper noise on the right.
  \label{fig:other-noise}
}

Because our recorded image $\bvec$ is a good approximation of  $\bvec^{\true}$, we might naively expect to find a good approximation of $\xvec^{\true}$ by solving the linear system of equations $\bA\xvec = \bvec$.  
However, its solution, which we call $\xvec^{LS}$, turns out to be very different from the original image $\xvec^{\true}$ as is seen in Figure~\ref{fig:naive-solution}.  As we will soon see, this behavior
results from the fact that deblurring is a discrete linear ill-posed problem.

\begin{myfigure}{
    \centering
    \includegraphics[scale=0.45]{Figures/64x64Case/x_true.png}
    \includegraphics[scale=0.45]{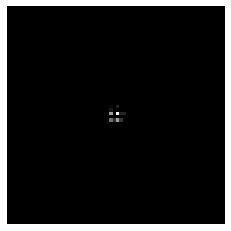}
  }
\end{myfigure}
\caption{ On the left we see the original image $\xvec^{\true}$ while
  the right shows $\xvec^{LS}$, the solution to the equation $\bA\xvec = \bvec$, where $\bvec$ is the noisy recorded image.
  \label{fig:naive-solution}
}

We are now faced with two questions: how can we reconstruct the original image $\xvec^{\true}$ more faithfully and how can we do it in a computationally efficient way.  For instance, today's typical phone photos have around 10 million pixels, which means that the number of entries in the blurring matrix $\bA$ is about 100 trillion.  Working with a matrix of that size will require some careful thought.

\section{Discrete Linear Ill-Posed Problems}

The singular value decomposition (SVD) of the matrix $\bA$ offers insight into 
why the naively reconstructed image $\bx^{LS}$ differs so greatly from
the original $\bx^{\true}$.  Consider both vectors $\bx^{\true}$ and $\bb$ of size $pq$ and define $m=pq$. Now suppose that $\bA\in\real^{m\times m}$ has
full rank and its SVD is given by
$$
\bA = \bU\bSigma \bV^T = \sum_{\ell=1}^m \sigma_\ell\bu_\ell\bv_\ell^T,
$$
where
$\bU=(\bu_1,\ldots,\bu_m)\in\real^{m\times m}$ and
$\bV=(\bv_1,\ldots,\bv_m)\in\real^{m\times m}$ are matrices having
orthonormal columns so that $\bU^T\bU = \bV^T\bV=\bI$, and
$$
\bSigma = {\rm diag}(\sigma_1,\ldots,\sigma_m),
\sigma_1\geq\sigma_2\geq\ldots\geq\sigma_m > 0.
$$
The scalars $\sigma_\ell$ are the {\em singular values} of $\bA$ and the
vector $\bu_\ell$ and $\bv_\ell$ are the {\em left} and {\em right singular
  vectors} of $\bA$, respectively. The singular value decomposition provides orthonormal bases defined by
the columns of $\bU$ and $\bV$ so that $\bA$ acts as scalar multiplication
by the singular values: $\bA\bv_\ell = \sigma_\ell\bu_\ell$.  Since the singular
values form a non-increasing sequence, the decomposition concentrates
the most important data in the beginning singular vectors, an observation
that is key to our work in reconstructing $\bx^{\true}$.

Discrete linear ill-posed problems are characterized by three
properties, which are satisfied by the deblurring problem.  First, the
singular values $\sigma_\ell$ decrease to zero without a large gap separating a group of large singular values from
a group of smaller ones.  The plot of the singular values in Figure \ref{fig:sigma-i} shows that the difference between the largest and smallest singular values is about thirteen orders of magnitude, which indicates that $\bA$ is highly ill-conditioned.

\begin{myfigure}{
    \centering
    \includegraphics[scale=0.4]{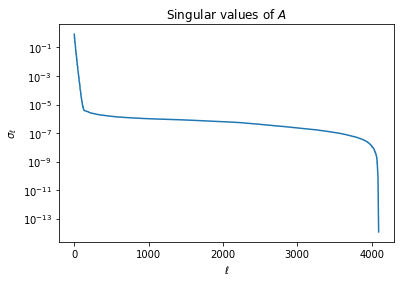}}
\figcaption{The singular values $\sigma_\ell$ of $\bA$.}
 \label{fig:sigma-i}
\end{myfigure}
  
Second, the singular vectors $\bu_\ell$ and $\bv_\ell$ become more and more
oscillatory as $\ell$ increases.  Figure \ref{fig:v-i} shows images $\bV_\ell$ representing eight right singular vectors $\bv_\ell$ ($\bv_\ell={\rm vec}(\bV_\ell)$) of the blurring matrix constructed above and demonstrates how the frequency of the oscillations increases as $\ell$
increases.  Since the blurring matrix $\bA$ spreads out any peaks in an image, it tends to dampen high frequency
oscillations. Therefore, a right singular vector $\bv_i$ representing a high frequency will correspond to a small singular value $\sigma_\ell$ since $\bA\bv_\ell = \sigma_\ell\bu_\ell$.

\begin{myfigure}{    
\centering
    \includegraphics[scale=0.45]{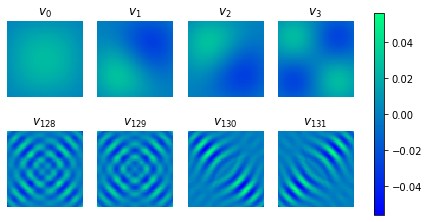}}
\end{myfigure}
\figcaption{Eight right singular vectors $\bv_\ell$. \label{fig:v-i}
}

The third property of discrete linear ill-posed problems is known as
the {\em discrete Picard condition}, which says that the coefficients
of $\bb^{\true}$ expressed in the left singular basis, $|\bu_\ell^T\bb^{\true}|$,
decay as fast as
the singular values $\sigma_\ell$.  This is illustrated on the left of
Figure \ref{fig:picard}, which shows both $|\bu_\ell^T\bb^{\true}|$ and the
singular values $\sigma_\ell$.  Since
$\bu_\ell^T\bb^{\true}=\sigma_\ell\bv_\ell^T\bx^{\true}$, the discrete Picard condition
holds when the original image $\bx^{\true}$ is not dominated by high
frequency contributions $\bv_\ell^T\bx^{\true}$ for large $\ell$.  This is a reasonable assumption with most digital
photographs. 

\begin{myfigure}{
  \centering
    \includegraphics[scale=0.25]{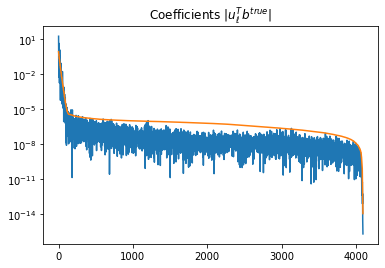}
    \includegraphics[scale=0.25]{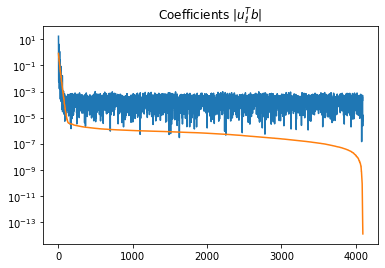}
  }
\end{myfigure}
\figcaption{ The coefficients $|\bu_\ell^T\bb^{\true}|$ are seen on the left while the coefficients $|\bu_\ell^T\bb|$ are on the right.
  \label{fig:picard}
}

By contrast, the coefficients $|\bu_\ell^T\bb| = |\bu_\ell^T\bb^{\true} + \bu_\ell^T\evec|$
appear on the right of Figure \ref{fig:picard}.  
The discrete Picard condition means that the coefficients $|\uvec_\ell^T\bvec^{\true}|$ decrease to zero.  However, the coefficients of the added white noise $|\uvec_\ell^T\evec|$ remain relatively constant for large $i$.  At some point the noise contributed by $\evec$ overwhelms the information contained in $\bvec^{\true}$.

Writing both vectors $\bx^{LS}$ and $\bx^{\true}$ as a
linear combination of the right singular vectors $\bv_\ell$, we see that
$$
\bx^{\true} = \sum_{\ell=1}^{m} \frac{\bu_\ell^T \bb^{\true}}{\sigma_\ell} \bv_\ell
$$
and
$$
\begin{aligned}
\bx^{LS}&= \sum_{\ell=1}^{m} \frac{\bu_\ell^T \bb}{\sigma_\ell} \bv_\ell =\sum_{\ell=1}^{m}
\frac{\bu_\ell^T (\bb^{\true}+\bee)}{\sigma_\ell} \bv_\ell \\ 
&= \sum_{\ell=1}^{m} \frac{\bu_\ell^T \bb^{\true}}{\sigma_\ell} \bv_\ell +
\sum_{\ell=1}^{m} \frac{\bu_\ell^T \bee}{\sigma_\ell} \bv_\ell \\
& = \bx^{\true} + \sum_{\ell=1}^{m} \frac{\bu_\ell^T \bee}{\sigma_\ell} \bv_\ell.
\end{aligned}
$$
Because the singular values approach
0, the coefficients $|\bu_\ell^T\bee|/\sigma_\ell$ grow extremely large, as seen in Figure \ref{fig:VTx-ls}.  Therefore, $\xvec^{LS}$ includes a huge contribution from high-frequency right singular vectors $\vvec_\ell$,
which means that $\bx^{LS}$ is
very oscillatory and not at all related to the original image
$\bx^{\true}$ that we seek to reconstruct.

\begin{myfigure}{
  \centering
    \includegraphics[scale=0.4]{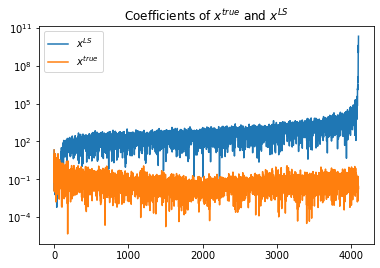}
  }
\end{myfigure}
\figcaption{
  The coefficients of $\xvec^{\true}$ and $\xvec^{LS}$.
  \label{fig:VTx-ls}
}

We also note here that the size of the coefficients shown in Figure \ref{fig:VTx-ls} cause the norm $\|\xvec^{LS}\|_2$ to be extremely large.  As we will see shortly, we will consider this norm as a measure of the amount of noise in the reconstructed image.

\section{Regularization}

As an alternative to accepting $\xvec^{LS}$ as our reconstructed image, we can compute approximations to $\bx^{\true}$
by filtering out the noise $\bee$ while retaining as much information as
possible from the measured data $\bb$.  This process is known as {\em regularization}, and there are several possible approaches we can follow. 

A first natural idea is to introduce a set of {\em filtering factors} $\phi_\ell$ on the SVD expansion and construct a
regularized solution as
$$
\bx^{\reg} = \sum_{\ell=1}^{m} \phi_\ell\frac{\bu_\ell^T \bb}{\sigma_\ell} \bv_\ell,
$$
with the filter factors being $\phi_\ell \approx 0$ for large values of $\ell$ when the noise dominates and $\phi_\ell \approx 1$ for small values of $\ell$, which are the terms in the expansion where the components of both $\bb$ and $\bb^{\rm true}$ are the closest.
  
One option is to define $\phi_\ell=1$ for $\ell$ smaller than some cut-off $k<m$ and  $\phi_\ell=0$ otherwise.
That is, we can simply truncate the expansion of $\bx^{LS}$ in terms of
right singular vectors in an attempt to minimize the contribution from
the terms $|\bu_\ell^T\evec|/\sigma_\ell$ for large $\ell$. Then, the obtained regularized solution would be
$$
\bx^{\reg} = \sum_{\ell=1}^{k} \frac{\bu_\ell^T \bb}{\sigma_\ell} \bv_\ell.
$$
This solution is known as the truncated SVD (TSVD). 

Remember, however, that the singular values in a discrete ill-posed problem approach 0 without there being a gap that would form a natural cut-off point. Instead, Tikhonov regularization chooses the filtering factors
$$
\phi_\ell = \frac{\sigma_\ell^2}{\sigma_\ell^2 + \lambda^2}
$$
for some parameter $\lambda$ whose choice will be discussed later.  For now,
notice that $\phi_\ell \approx 1$ when $\sigma_\ell\gg \lambda$ and $\phi_\ell
\approx 0$ when $\sigma_\ell\ll \lambda$.  This has the effect of truncating
the singular vector expansion at the point where the singular values pass through
$\lambda$ but doing so more smoothly.

This leads to the regularized solution
$$
\bx^{\reg} = \sum_{\ell=1}^{m} \frac{\sigma_\ell}{\sigma_\ell^2 + \lambda^2}(\bu_\ell^T \bb) \bv_\ell,
$$
which may also be rewritten as
$$
\bx^{\reg} = (\bA^T\bA+\lambda^2\bI)^{-1}\bA^T\bb.
$$
This demonstrates that $\bx^{\reg}$ solves the least squares problem
\begin{equation}\label{eq:Tikhonov}\bx^{\reg}=\argmin_{\bx} \{\|\bA\bx -\bb\|_2^2 + \lambda^2\|\bx\|_2^2\}.\end{equation}

This is a helpful reformulation of Tikhonov regularization.  First, writing $\xvec^{\reg}$ as the solution of a least squares problem provides us with efficient computational alternatives to finding the SVD of $\bA$.  Moreover, the minimization problem provides insight into choosing the optimal value of the regularization parameter $\lambda$, as we will soon see.

By the way, this formulation of Tikhonov regularization shows its
connection to ridge regression, a data science technique for tuning a
linear regression model in the presence of multicollinearity to improve its predictive accuracy.

Let us investigate the meaning of (\ref{eq:Tikhonov}).  Notice that
\begin{equation}
    \label{eq:true-residual}
\|\bA\xvec^{\true} - \bvec\|_2^2 = \|\bvec^{\true} - \bvec\|_2^2 = \|\evec\|_2^2,
\end{equation}
which is relatively small.  We therefore consider the residual $\|\bA\xvec-\bvec\|_2^2$ as a measure of how far away we are from the original image $\xvec^{\true}$.  Remember that $|\uvec_\ell^T\evec|/\sigma_\ell$, the contributions to $\xvec^{LS}$ from the added noise, cause $\|\xvec^{LS}\|_2$ to be very large.  Consequently, we think of the second term in (\ref{eq:Tikhonov}) as measuring the amount of noise in $\xvec$.

The regularization parameter $\lambda$ allows us to balance these two sources of error.  For instance, when $\lambda$ is small, the regularized solution $\xvec^{\reg}$, which is the minimum of (\ref{eq:Tikhonov}), will have a small residual $\|\bA\xvec^{\reg}-\bb\|_2$ at the expense of a large norm $\|\xvec^{\reg}\|_2$.  In other words, we tolerate a noisy regularized solution in exchange for a small residual.  On the other hand, if $\lambda$ is large, the regularized solution will have a relatively large residual in exchange for filtering out a lot of the noise.  

\section{Choosing the regularizing parameter $\lambda$}
When applying Tikhonov regularization for solving discrete ill-posed inverse problems, a question of high interest is:
How can we define the best regularization parameter
$\lambda$? A common and well known technique is
to create a log-log plot of the residuals $\|\bA\bx^{\reg}-\bb\|_2$ and the norm
$\|\bx^{\reg}\|_2$ as we vary $\lambda$. This plot, as shown in Figure
\ref{fig:L-curve}, is usually called an $L$-curve due to its
characteristic shape~\cite{hansen1993use}.

As mentioned earlier, small values of $\lambda$ lead to noisy regularized solutions while large values of $\lambda$ produce large residuals.  This means that we move across the $L$-curve from the upper left to the lower right as we increase $\lambda$.  

\begin{myfigure}{
  \centering
  \includegraphics[scale=0.45]{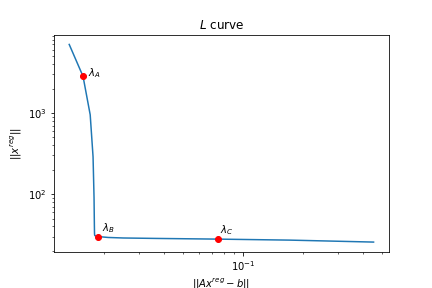}
}
\end{myfigure}
\figcaption{
  The $L$-curve in our sample deblurring problem.  The three indicated points correspond to values of $\lambda = \lambda_A, \lambda_B$, and $\lambda_C$.
  \label{fig:L-curve}
}

Since the filtering factors satisfy $\phi_i\approx 1$ when $\sigma\gg\lambda$ and $\phi_\ell\approx0$ when $\sigma_\ell\ll\lambda$, we will
view the point where $\sigma_\ell = \lambda$ as an indication of where we begin filtering.  
Let us first consider the value $\lambda = \lambda_A$, which produces the point on the $L$-curve indicated in Figure \ref{fig:L-curve}.  The resulting coefficients $|\uvec_\ell^T\bvec|$ and $\phi_\ell|\uvec_\ell^T\bvec|$ are shown in Figure \ref{fig:filtering-A}.  Filtering begins roughly where the plot of singular values crosses the horizontal line indicating the value of $\lambda_A$.  

\begin{myfigure}{
  \centering
  \includegraphics[scale=0.5]{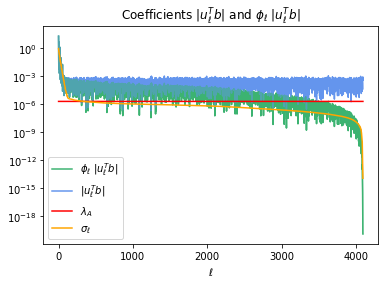}
}
\end{myfigure}
\figcaption{
  The choice $\lambda=\lambda_{A}$ leads to under-smoothing.
  \label{fig:filtering-A}
}

While we have filtered out some of the noise, it appears that there is still a considerable amount of noise present.  This is reflected by the position of the corresponding point on the $L$-curve since the norm $\|\xvec^{\reg}\|_2$ is relatively large.  This choice of $\lambda$ is too low, and we say that the regularized solution is under-smoothed.

Alternatively, let us consider the regularized solution constructed with $\lambda = \lambda_C$ as indicated on Figure \ref{fig:L-curve}.  This leads to the coefficients $|\uvec_\ell^T\bvec|$ and $\phi_\ell|\uvec_\ell^T\bvec|$ shown in Figure \ref{fig:filtering-C}.

\begin{myfigure}{
  \centering
  \includegraphics[scale=0.5]{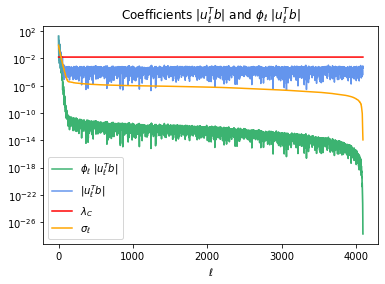}
}
\end{myfigure}
\figcaption{
  The choice $\lambda=\lambda_C$ leads to over-smoothing.
  \label{fig:filtering-C}
}

In this case, we begin filtering too soon so that, while we have removed the noise, we have also discarded some of the information present in $\bvec$, which is reflected in the relatively large residual $\|\bA\xvec^{\reg}-\bvec\|_2$. This choice of $\lambda$ is too large, and we say the regularized solution is over-smoothed.

Finally, considering the case where $\lambda = \lambda_B$ gives the regularized solution that appears at the sharp bend of the $L$-curve in Figure \ref{fig:L-curve}.  The resulting coefficients shown in Figure \ref{fig:filtering-B} inspire confidence that this is a good choice for $\lambda$.

\begin{myfigure}{
  \centering
  \includegraphics[scale=0.5]{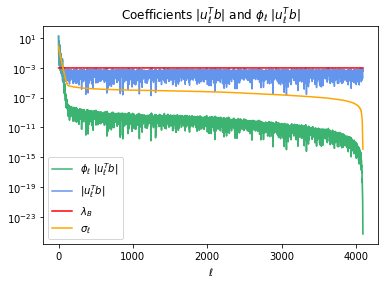}
}
\end{myfigure}
\figcaption{
  The choice $\lambda=\lambda_B$ gives the optimal amount of smoothing.
  \label{fig:filtering-B}
}

In this case, decreasing $\lambda$ causes us to move upward on the $L$-curve;  we are adding noise without improving the residual.  Increasing $\lambda$ causes us to move right on the $L$-curve;  we are losing information as the residual increases without removing any more noise.  Therefore, $\lambda = \lambda_B$ is our optimal value.

Figure \ref{fig:x-reg} shows the image $\bx^{\reg}$ obtained using this
optimal parameter $\lambda=\lambda_B$.  While it is not a perfect reconstruction of the
original image $\bx^{\true}$, it is a significant improvement over the
recorded image $\bb$.  Accurately reproducing the sharp boundaries that occur in $\xvec^{\true}$ requires high-frequency contributions that we have necessarily filtered out.  

\begin{myfigure}{
  \centering
  \includegraphics[scale=0.45]{Figures/64x64Case/x_true.png}
  \includegraphics[scale=0.45]{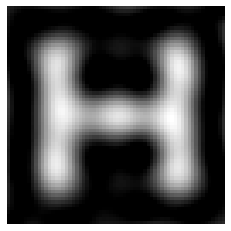}
}
\end{myfigure}
\figcaption{
  The reconstructed image $\bx^{\reg}$ using the regularization
  parameter determined by the $L$-curve is seen on the right, along
  with the original image $\bx^{\true}$ on the left for comparison.
  \label{fig:x-reg}
}

The $L$-curve furnishes a practical way to identify the optimal regularizing parameter as there are techniques that allow us to identify the point of maximal curvature by computing the regularized solution for just a few choices of $\lambda$.  However, this technique should not be applied uncritically as there are cases in which the optimal regularized solution does not converge to the true image as the added error approaches zero.

An alternative technique, known as the Discrepancy Principle~\cite{engl1987discrepancy}, relies on an estimate of the size of the error $\|\evec\|_2$.  Remember from (\ref{eq:true-residual}) that we have $\|\bA\xvec^{\true} - \bvec\|_2 = \|\evec\|_2$.  Moreover, the SVD description of $\xvec^{\reg}$ provides a straightforward explanation for why the residual $\|\bA\xvec^{\reg} - \bvec\|_2$ is an increasing function of $\lambda$.  If we know $\|\evec\|_2$, we simply choose the optimal $\lambda$ to be the one where $\|\bA\xvec^{\reg} - \bvec\|_2 = \|\evec\|_2$.

Other well-known methods for choosing the parameter $\lambda$ include the Generalized Cross Validation (GCV)~\cite{golub1979generalized}, that chooses $\lambda$ to maximize the accuracy with which we can predict the value of a pixel that has been omitted, the unbiased predictive risk estimator (UPRE)~\cite{ReVaAr17}, and more recently, methods based on learning when training data is available~\cite{chung2021efficient, chung2017learning}.

\section{Other Regularization Techniques}
Looking at the variational formulation \eqref{eq:Tikhonov} of Tikhonov regularization, it is easy to see how it can be extended to define other regularization methods by, for example, using different regularization terms.

\subsection{General Tikhonov Regularization}
The Tikhonov regularization formulation \eqref{eq:Tikhonov} can be generalized to 
\begin{equation}\label{eq: GTik}
    \min_{\bx} \{\|\bA\bx - \bb\|^2_2 + \lambda\|\bL\bx\|^2_2\},
\end{equation}
by incorporating a matrix $\bL$, which is called the regularization matrix. Its choice is problem dependent and can significantly affect the quality of the reconstructed solution. Several choices of the regularization matrix involve discretization of the derivative operators or framelet and wavelet transformations depending on the application. The only requirement on $\bL$ is that it should satisfy
\[
\mathcal{N}(\bA)\cap\mathcal{N}(\bL)=\{0\},
\]
where $\mathcal{N}(\bM)$ denotes the null space of the matrix $\bM$. The general Tikhonov 
minimization problem \eqref{eq: GTik} has the unique solution 
\[
\bx_\lambda=(\bA^T\bA+\lambda \bL^T\bL)^{-1}\bA^T\bb
\]
for any $\lambda>0$. 

\begin{myfigure}{
  \centering
  \includegraphics[scale=0.45]{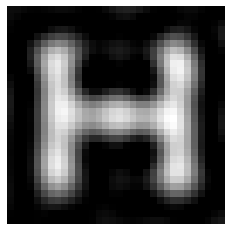}
  \includegraphics[scale=0.45]{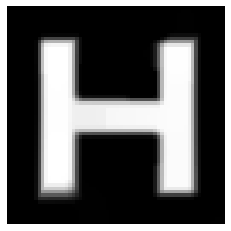}
}
\end{myfigure}
\figcaption{
  The reconstructed images $\bx^{\reg}$ using the optimal regularization
  parameter and the discretization of the first derivative operator with the 2-norm regularization on the left, and TV regularization on the right. 
  \label{fig:xgTik}
}

In Figure \ref{fig:xgTik}, we reconstruct the image by applying the discretization of the two-dimensional first derivative operator for zero boundary conditions, that is, the matrix $\bL$ takes the form of 

\begin{equation}\bL=\left( \begin{matrix} \bI&\otimes&\bL_1\\ \bL_1&\otimes& \bI\end{matrix}\right) \mbox{with }
\bL_1=\left( \begin{smallmatrix}-1 & 1&&&\\
&-1&1&&\\
&&\ddots&\ddots&\\
&&&-1&1 \end{smallmatrix}\right)\label{eq:Lmatrix}\end{equation}
where $\otimes$ is the \emph{Kronecker product} \cite{GVL12} defined, for matrices $\bB$ and $\bC$, by
$$ \bB\otimes \bC=\left( \begin{matrix} b_{11} \bC & b_{12} \bC & \dots & b_{1m}\bC\\ b_{21} \bC & b_{22} \bC & \dots &  b_{2m} \bC\\ \vdots&\vdots&& \vdots\\  b_{m1} \bC&  b_{m2} \bC & \dots & b_{mm} \bC
\end{matrix}\right).$$

\subsection{Total Variation Regularization}
In many applications, the image to be reconstructed is known to be piece-wise constant with regular and sharp edges, like the one we are using as an example in this article. Total variation (TV) regularization is a popular choice that allows the solution to preserve edges. Such regularization can be formulated as $$\min_\bx\{\|\bA\bx - \bb\|^2_2 + \lambda\|\bL\bx\|_1\},$$ where $\bL$ is the TV operator, which once discretized, it is the same as \eqref{eq:Lmatrix}. Looking at Figure \ref{fig:xgTik}, we can see that even though we are using the same operator $\bL$, the norms used in the regularization terms are different and that makes a huge difference. But there is a higher cost to finding the TV solution, due to the fact that the minimization functional is not differentiable. Still, there are many algorithms to find its minimum. Here, we apply the Iteratively Reweighted Least Squares (IRLS), that solves a sequence of general form Tikhonov problems. So, instead of solving only one Tikhonov problem, we solve many. 

The TV approach is also very commonly used in compressed sensing where the signal to be reconstructed is sparse in its original domain or in some transformed domain~\cite{candes2006stable,lustig2008compressed}. Recently it has been used in the context of regularizing large-scale dynamic inverse problems \cite{pasha2021efficient} as well as in learning when training data is available \cite{antil2020bilevel}.

\section{Matrix Structures}\label{sec: matrix_structure}
\subsection{BTTB Structure}
Because of the large-scale problems, it is useful to consider the structure of the matrix $\bA$. For instance, when considering a spatially invariant blur and assuming that the image has zero boundary conditions, the matrix $\bA$ has Block-Toeplitz-Toeplitz-Block (BTTB) structure \cite{hansen1994regularization}, that is,
a $p$ by $p$ block-Toeplitz matrix with each block being a $p$ by $p$ Toeplitz matrix,
$$\bA=\left( \begin{matrix}
\bA_0 & \bA_{-1} & \bA_{-2}& \dots & \bA_{-(p-1)}\\
\bA_1 & \bA_0&\bA_{-1}& \dots&\bA_{-(p-2)}\\
\bA_2&\bA_1&\bA_0& \dots & \bA_{-(p-3)}\\ \vdots & \vdots
& \vdots & \ddots &\vdots \\
\bA_{p-1}&\bA_{p-2}&\bA_{p-3}& \dots &\bA_0
\end{matrix}\right),$$
where for $\ell=1,\dots, p$ we have
$$\bA_\ell=\left( \begin{matrix}
a^\ell_0 & a^\ell_{-1} & a^\ell_{-2}& \dots & a^\ell_{-(p-1)}\\
a^\ell_1 & a^\ell_0&a^\ell_{-1}& \dots&a^\ell_{-(p-2)}\\
a^\ell_2&a^\ell_1&a^\ell_0& \dots & a^\ell_{-(p-3)}\\ \vdots & \vdots
& \vdots & \ddots &\vdots \\
a^\ell_{p-1}&a^\ell_{p-2}&a^\ell_{p-3}& \dots &a^\ell_0
\end{matrix}\right).$$

Notice that to generate these matrices, we only need the
first row and column of each matrix $\bA_\ell$, what is called the
Toeplitz vector $\baa^\ell= (a^\ell_{-(p-1)}, \dots, a^\ell_{p-1})$.

\subsection{BCCB Structure}
Another feasible structure arises when still considering a
spatially invariant blur but assuming that the image has
periodic boundary conditions. In this case, the matrix $\bA$
has a Block-Circulant-Circulant-Block (BCCB) structure,
that is, a $p$ by $p$ block-circulant matrix with each block
being a $p$ by $p$ circulant matrix,

$$\bA=\left( \begin{matrix}
\bA_1 & \bA_2 & \bA_3& \dots & \bA_p\\
\bA_p & \bA_1&\bA_2& \dots&\bA_{p-1}\\
\bA_{p-1}&\bA_p&\bA_1& \dots & \bA_{p-2}\\ \vdots & \vdots
& \vdots & \ddots &\vdots \\
\bA_2&\bA_3&\bA_4& \dots &\bA_1
\end{matrix}\right),$$
where for $\ell=1,\dots, p$
$$\bA_\ell=\left( \begin{matrix}
a^\ell_1 & a^\ell_2 & a^\ell_3& \dots & a^\ell_p\\
a^\ell_p & a^\ell_1&a^\ell_2& \dots&a^\ell_{p-1}\\
a^\ell_{p-1}&a^\ell_p&a^\ell_1& \dots & a^\ell_{p-2}\\ \vdots & \vdots
& \vdots & \ddots &\vdots \\
a^\ell_2&a^\ell_3&a^\ell_4& \dots &a^\ell_1
\end{matrix}\right).$$

Notice that to generate these matrices, we only need the
first row of each matrix $\bA_\ell$, ${\bf a}^\ell = (a^\ell_1,\dots, a^\ell_
p)$, so we
do not need to store every entry of the matrix $\bA$, only
the vectors that define the matrices $\bA_\ell$. We might never need to
build the matrix $\bA$ either to get a regularized solution. Furthermore, BCCB
matrices are diagonalizable by the two-dimensional Discrete Fourier Transform (DFT) with eigenvalues
given by the Fourier transform of the first column of
$\bA$. For instance, the Tikhonov solution can be written as
$${\bf X}_\lambda= \mbox{IDFT} \left(\left(\frac{\mbox{conj}(\hat{\bf a})}{|\hat{\bf a}|^2+\lambda {\bf 1}}\right)\odot \mbox{DFT}(\bB)\right),$$
where $\bX_\lambda$ and $\bB$ are the matrices representing $\bx_\lambda$ and $\bb$, respectively (i.e., $\bx_\lambda={\rm vec} (\bX_\lambda)$ and $\bb={\rm vec}(\bB)$), $\hat{\bf a}=m \mbox{DFT}({\bf a}_s)$, conj($\cdot$) denotes the componentwise complex conjugate, $|\hat{\bf a}|^2$ is the $m\times m$ matrix whose entries are the squared magnitudes of the complex entries of $\hat{\bf a}$, ${\bf 1}$ is a matrix of all ones, and $\odot$ denotes component-wise multiplication between matrices. Notice that IDFT stands for Inverse Discrete Fourier Transform.
For more details, we direct the reader to \cite[Chapter~5]{vogel2002computational}.

\subsection{Separable Blur Operator}
Another special case is when the blurring operator $\bA$ is
separable, which means that the blurring can be separated
in the horizontal and vertical directions as illustrated in Equation \eqref{eq:Gaussian-blur} that is, the matrix
$\bA$ can be written as
$$\bA = \bA_r \otimes \bA_c.$$

If we have the SVDs $\bA_r=\bU_r\bSigma_r \bV^T_r$ and $\bA_c=\bU_c\bSigma_c \bV^T_c,$ we basically have the SVD of $\bA$ by writing
\begin{align}
\bA &= \bA_r \otimes \bA_c \\
&=(\bU_r\bSigma_r \bV_r^T)\otimes (\bU_c\bSigma_c \bV_c^T) \\
&=(\bU_r\otimes \bU_c)(\bSigma_r\otimes \bSigma_c)(\bV_r\otimes \bV_c)^T,
\end{align}
except that the elements in the diagonal matrix $\bSigma_r\otimes \bSigma_c$ might not be in the decreasing order, and therefore some reordering might be needed. 
So, the Tikhonov solution can be written as
$${\bf X}_\lambda={\bf V}_c\left(\left(\frac{\bd_c\bd_r^T}{(\bd_c\bd_r^T)^2+\lambda {\bf 1}}\right)\odot(\bU_c^T\bB\bU_r)\right) \bV_r^T,$$
where ${\bf 1}$ is a matrix of all ones, and $\bd_r$ and $\bd_c$ are the diagonals of $\bSigma_r$ and $\bSigma_c$, respectively.

Notice that if $\bA_r$ and $\bA_c$ are Toeplitz matrices, then $\bA_r\otimes \bA_c$ has BTTB structure, and if $\bA_r$ and $\bA_c$ are circulant
matrices, then $\bA_r\otimes \bA_c$ has BCCB structure.

\section{Multilevel Methods}
Because dealing with images of large size is difficult, researchers keep working on finding ways to solve these large-scale inverse problems efficiently. One possible way is by developing multilevel methods. The main idea of a multilevel method is to define a sequence of systems of equations decreasing in size,
$$\bA^{(n)}\bx^{(n)} = \bb^{(n)}, \quad  0\leq n\leq L,$$
where the superscript $n$ denotes the $n$-th level ($n$ = 0 being
the original system), and to get an approximate solution, a correction, 
or some information at each level where the computational
cost would be smaller than solving the original
system.  At each level, the right-hand side is defined by
$$\bb^{(n+1)} = \bR^{(n)}\bb^{(n)}$$ and the matrix by
$$\bA^{(n+1)} = \bR^{(n)}\bA^{(n)}\bP^{(n)},$$
where $\bR^{(n)}$ is called the restriction operator and $\bP^{(n)}$ is the interpolation
operator.

There are many ways of defining and using this
sequence of systems of equations to solve many different
mathematical problems. To learn more about multilevel methods for image deblurring applications,
we recommend reference \cite{chan2010multilevel, donatelli2006regularizing, espanol2010multilevel, morigi2008cascadic}

\subsection{Wavelet-based approach} Here, we consider
the use of wavelet transforms as restriction
and interpolation operators.  In particular, we will work with the Haar Wavelet Transform (HWT) because, as we will see, it keeps the structures of the matrices involved. The one-dimensional HWT is the
$p\times p$ orthonormal matrix $\bW$ defined by
$$\bW=\frac{1}{\sqrt{2}}\left( \begin{smallmatrix}
1 & 1 & 0 & 0 & \dots & \dots & 0 & 0\\
0 & 0 & 1 & 1 & \dots & \dots & 0 & 0\\
\vdots & \vdots & \vdots & \vdots & \ddots & \ddots & \vdots & \vdots\\
0 & 0 & 0 & 0 & \dots & \dots & 1 & 1\\
1 & -1 & 0 & 0 & \dots & \dots & 0 & 0\\
0 & 0 & 1 & -1 & \dots & \dots & 0 & 0\\
\vdots & \vdots & \vdots & \vdots & \ddots & \ddots & \vdots & \vdots\\
0 & 0 & 0 & 0 & \dots & \dots & 1 & -1\\
\end{smallmatrix}\right) = \left( \begin{matrix}\bW_1\\\bW_2\end{matrix}\right),$$
where $\bW_1$ and $\bW_2$ have dimension $p/2 \times p$. 
The two-dimensional HWT can be defined in terms of the one-dimensional HWT by ${\bW_{\rm 2D}} = \bW \otimes \bW \in \mathbb{R}^{m\times m}$ with $m=p^2$. So, we define the restriction operator by $\bR = \bW_1 \otimes \bW_1 \in  \mathbb{R}^{(m/4)\times m}$ and the interpolation operator by $\bP = \bR^T$.

To motivate the use of HWT as a restriction operator, on the left of Figure \ref{fig:comp}, we can see the coarser version of $\bx^{\rm true}$ defined by $\bx^{(1)}=\bR\bx^{\rm true}$.  This is a $32\times 32$ image that still shows basically the same information as the original $64\times 64$ image with the letter H. So, reconstructing a coarse version could be enough for some tasks. Another reason is that, we actually do this many times when \emph{compressing} images to save storage space in the computer. So imagine getting the image $\bb$, and before transmitting it, we compute its coarser version $\bb^{(1)}=\bR\bb$ and transfer just $\bb^{(1)}$ (see right of Figure \ref{fig:comp}). The image $\bb^{(1)}$ might have enough information to recover $\bx^{(1)}$.

\begin{myfigure}{
  \centering
  \includegraphics[scale=0.45]{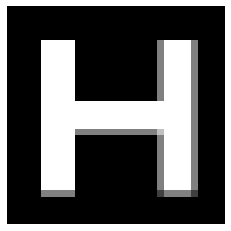}
  \includegraphics[scale=0.45]{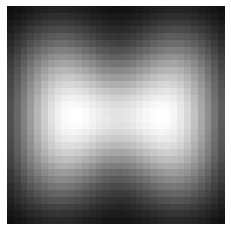}
}
\end{myfigure}
\figcaption{
  The compressed true image $\bx^{(1)}$ on the left and blurred and noisy image $\bb^{(1)}$ on the right. 
  \label{fig:comp}
}

The question now is what is the right blurring matrix to recover $\bx^{(1)}$ from $\bb^{(1)}$. Using our example, we will show that using $\bA^{(1)} = \bR\bA\bP$ does a good job. Figure \ref{fig:xMGM} shows the Tikhonov and TV solutions of the system $\bA^{(1)}\bx = \bb^{(1)}$. Notice that we are solving a system with a $1,024 \times 1,024$ matrix instead of the original system with a $4,096\times 4,096$ matrix.

In the following theorems, we want to show that the use of HWT keeps the nice structures of matrices mentioned before.

\begin{theorem} \cite{espanol2010multilevel} Let $\bT$ be a $p\times p$ matrix with Toeplitz
structure and with Toeplitz vector $\bt$, and $p = 2^s$. Then,
the $2^{s-1}\times 2^{s-1}$ matrix $\bT^1 = \bW_1\bT\bW^T_1$ is also Toeplitz
with Toeplitz vector $\bt_1 = \tilde \bt(1: 2: 2p - 3)$, where $\tilde \bt = \tilde \bT \bt$
with $\tilde \bT$ being the Toeplitz matrix with Toeplitz vector with
all zeros and $\bt_0 = \bt_{-2} = 1/2$ and $\bt_{-1} = 1$.
\end{theorem}
\begin{myfigure}{
  \centering
  \includegraphics[scale=0.45]{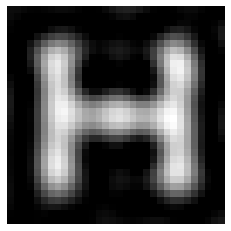}
  \includegraphics[scale=0.45]{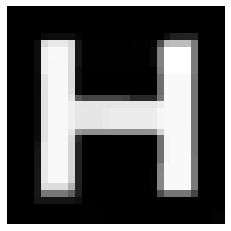}
}
\end{myfigure}
\figcaption{
  The reconstructed coarse images $\bx^{\reg}$ of size $32\times 32$, obtained using the optimal regularization
  parameter and the discretization of the first derivative operator with the 2-norm regularization on the left, and TV regularization on the right. 
  \label{fig:xMGM}
}

A similar result can be shown for circulant matrices.
\begin{corollary} Let $\bC$ be a $p\times p$ circulant matrix and
$p = 2^s$. Then, the $2^{s-1} \times 2^{s-1}$ matrix $\bC_1 = \bW_1\bC \bW^T_1$
is also circulant.
\end{corollary}

Let us consider the case when $\bA = \bA_r \otimes \bA_c \in \mathbb{R}^{m\times m}$. Then,
by property of the Kronecker product we have that
\begin{align}\bA^{(1)} &= (\bW_1 \otimes \bW_1) (\bA_r \otimes \bA_c) (\bW_1 \otimes  \bW_1)^T \nonumber\\
&= (\bW_1A_r\bW^T_1)\otimes (\bW_1\bA_c\bW^T_1)\nonumber \\
&= \bA^{(1)}_r\otimes \bA^{(1)}_c. \nonumber\end{align}

The matrix $\bA^{(1)}\in \mathbb{R}^{(m/4)\times(m/4)}$ is separable. Furthermore,
by Theorem 1, if $\bA_r$ and $\bA_c$ are Toeplitz matrices, then $\bA^{(1)}_r$
and $\bA^{(1)}_c$ are too. Applying this same argument again
and again, we obtain that $\bA^{(n)}$ is separable by two Toeplitz
matrices and therefore BTTB for all levels $n = 0,\dots, L$.
Similarly, by Corollary 1, $\bA^{(n)}$ is separable by two circular
matrices and therefore BCCB for all levels $n = 0,\dots, L$.
Therefore, the initial structure of the matrix is inherited to
all the levels.

For the case when we have BCCB structures, we can
solve the corresponding Tikhonov systems at all levels
using Fourier-based methods. If we are dealing with
separable matrices with Toeplitz structure, we could go
down several levels until we can compute the SVD of $\bA^{(n)}_r=\bU^{(n)}_r\bSigma^{(n)}_r(\bV^{(n)}_r)^T$ and $\bA^{(n)}_c=\bU^{(n)}_c\bSigma^{(n)}_c(\bV^{(n)}_c)^T,$ and use that
\begin{align}
\bA^{(n)} &= \bA^{(n)}_r \otimes \bA^{(n)}_c \nonumber \\
&=(\bU^{(n)}_r\bSigma^{(n)}_r(\bV^{(n)}_r)^T)\otimes (\bU^{(n)}_c\bSigma^{(n)}_c(\bV^{(n)}_c)^T) \nonumber\\
&=(\bU^{(n)}_r\otimes \bU^{(n)}_c)(\bSigma^{(n)}_r\otimes \bSigma^{(n)}_c)(\bV^{(n)}_r\otimes \bV^{(n)}_c)^T. \nonumber
\end{align}
This gives us basically the SVD of $\bA^{(n)}$, except that the elements in the diagonal matrix $\bSigma^{(n)}_r\otimes \bSigma^{(n)}_c$ might not be in the decreasing order, and therefore some reordering would be needed. 

There are also other efficient methods such as computing
the closest BCCB matrix and using it as preconditioner (see
\cite[Chapter~4]{vogel2002computational} or \cite{hansen2006deblurring} for more details). 

\section{Conclusions and outlook}

Our intention has been to introduce readers to the problem of image deblurring, the mathematical issues that arise, and a few techniques for addressing them.  As mentioned in the introduction, these techniques may be applied to a wide range of related problems.  Indeed, we outlined a number of alternative strategies, for example, in choosing appropriate boundary conditions or in selecting the best regularization parameter, as some strategies are better suited to a specific range of applications.  

We believe this subject is accessible to both undergraduate and graduate students and can serve as a good introduction to inverse problems, working with ill-conditioned linear operators, and large-scale computation.  The visual nature of the problem provides compelling motivation for students and allows the efficacy of various techniques to be easily assessed. The references ~\cite{ hansen2010discrete, hansen2006deblurring, vogel2002computational} contain excellent introductions to this subject.
Furthermore, the code to generate the figures appearing in this article is available at  \url{https://github.com/ipiasu/AMS_Notices_AEP}.

In addition, this field is an active area of research that aligns with the recent developments in machine learning and convolutional neural networks. Many classical techniques that have been traditionally used in image deblurring can serve as a tool to speed up the computations of the more recent methods that aim to learn models and parameters when training data are available or train a network to classify images. More particularly, we mentioned earlier in the manuscript that Tikhonov regularization is related to ridge regression, a fundamental tool in machine learning. The connection to machine learning goes much deeper. 
\end{multicols}

\section{Acknowledgments}
This article will appear on the AMS Notices. 
\bibliographystyle{plain}
\bibliography{biblio.bib}
\end{document}